\documentclass[twoside,11pt]{article}

\usepackage{jmlr2e}

\usepackage{amsmath}
\usepackage{amssymb}
\usepackage{delarray}
\usepackage{amsfonts}
\usepackage[english]{babel}
\usepackage{enumerate}

\newcommand \cA{{\cal A}}
\newcommand \cB{{\cal B}}
\newcommand \cC{{\cal C}}

\newcommand \cF{{\cal F}}

\newcommand \cX{{\cal X}}

\newcommand{\1}{{\rm 1}\kern-0.24em{\rm I}}

\jmlrheading{?}{2006}{??}{11/05}{??/??}{Guillaume Lecu\'e}

\newtheorem{theo}{Theorem}

\newtheorem{lem}{Lemma}

\newtheorem{defn}{Definition}

\ShortHeadings{Lower bounds and aggregation in density
estimation}{Lecu\'e} \firstpageno{1}

\begin{document}

\title{Lower bounds and aggregation in density
estimation}

\author{\name Guillaume Lecu\'e \email lecue@ccr.jussieu.fr \\
       \addr Laboratoire de Probabilit\'es et Mod\`eles Al\'eatoires\\
        Universit\'e Paris 6\\
        4 place Jussieu, BP 188\\
        75252 Paris, France}

\editor{G{\'a}bor Lugosi}

\maketitle

\begin{abstract} In this paper we prove the optimality of an aggregation
procedure. We prove lower bounds for aggregation of model
selection type of $M$ density estimators for the Kullback-Leiber
divergence (KL), the Hellinger's distance and the $L_1$-distance.
The lower bound, with respect to the KL distance, can be achieved
by the on-line type estimate suggested, among others, by
\citet{yang:00}. Combining these results, we state that $\log M/n$
is an optimal rate of aggregation in the sense of \citet{tsy:03},
where $n$ is the sample size.
\end{abstract}


\section{Introduction}
Let $(\cX,\cA)$ be a measurable space and $\nu$ be a
$\sigma$-finite measure on $(\cX,\cA)$. Let
$D_{n}=(X_{1},\ldots,X_{n})$ be a sample of $n$ i.i.d.
observations drawn from an unknown probability of density $f$ on
$\cX$ with respect to $\nu$. Consider the estimation of $f$ from
$D_n$.

Suppose that we have $M\geq 2$ different estimators
$\hat{f}_1,\ldots,\hat{f}_M$ of $f$. \citet{cat:97},
\citet{yang:00}, \citet{nem:00}, \citet{jn:00}, \citet{tsy:03},
\citet{catbook:04}  and \citet{rt:04} have studied the problem of
model selection type aggregation. It consists in construction of a
new estimator $\tilde{f}_n$ (called \emph{aggregate}) which is
approximatively at least as good as the best among
$\hat{f}_1,\ldots,\hat{f}_M$. In most of these papers, this
problem is solved by using a kind of cross-validation procedure.
Namely, the aggregation is based on splitting the sample in two
independent subsamples $D^{1}_{m}$ and $D^{2}_{l}$ of sizes $m$
and $l$ respectively, where $m \gg l$ and $m+l=n$. The size of the
first subsample has to be greater than the one of the second
because it is used for the true estimation, that is for the
construction of the $M$ estimators $\hat{f}_1,\ldots,\hat{f}_M$.
The second subsample is used for the adaptation step of the
procedure, that is for the construction of an aggregate $\tilde
f_{n}$, which has to mimic, in a certain sense, the behavior of
the best among the estimators $\hat{f}_i$. Thus, $ \tilde f_n $ is
measurable w.r.t. the whole sample $ D_n $ unlike the first
estimators $\hat{f}_1,\ldots,\hat{f}_M$.

One can suggest different aggregation procedures and the question
is how to look for an optimal one. A way to define optimality in
aggregation in a minimax sense for a regression problem is
suggested in \citet{tsy:03}. Based on the same principle we can
define optimality for density aggregation. In this paper we will
not consider the sample splitting and concentrate only on the
adaptation step, i.e. on the construction of aggregates (following
\citet{nem:00}, \citet{jn:00}, \citet{tsy:03}). Thus, the first
subsample is fixed and  instead of estimators
$\hat{f}_1,\ldots,\hat{f}_M$, we have fixed functions
$f_{1},\ldots,f_M $. Rather than working with a part of the
initial sample we will use, for notational simplicity, the whole
sample $D_n$ of size $n$ instead of a subsample $D_l^2$.

The aim of this paper is to prove the optimality, in the sense of
\citet{tsy:03}, of the aggregation method proposed by Yang, for
the estimation of a density on $(\mathbb{R}^d,\lambda)$ where
$\lambda$ is the Lebesgue measure on $\mathbb{R}^d$. This
procedure is a convex aggregation with weights which can be seen
in two different ways. Yang's point of view is to express these
weights in function of the likelihood of the model, namely
\begin{equation}\label{esti}
\tilde{f}_n(x)=\sum_{j=1}^M\tilde{w}_j^{(n)}f_j(x),\quad \forall
x\in\cX,
\end{equation}where the weights are
$\tilde{w}_j^{(n)}=(n+1)^{-1}\sum_{k=0}^nw_j^{(k)}$
 and \begin{equation}\label{weights}w_j^{(k)}=\frac{f_j(X_1)\ldots f_j(X_k)}{\sum_{l=1}^{M}f_l(X_1)\ldots
f_l(X_k)}, \ \forall k=1,\ldots,n \mbox{ and }
w_j^{(0)}=\frac{1}{M}.\end{equation} And the second point of view
is to write these weights as exponential ones, as used in
\citet{abb:97}, \citet{catbook:04}, \citet{har:02}, \citet{bn:05},
\citet{jntv:05} and \citet{lec:05}, for different statistical
models. Define the empirical Kullback loss $K_n(f)=-
(1/n)\sum_{i=1}^n \log{f(X_i)}$ (keeping only the term independent
of the underlying density to estimate) for all density $f$. We can
rewrite these weights as exponential weights:
$$w_j^{(k)}=\frac{\exp(-kK_k(f_j))}{\sum_{l=1}^M\exp(-kK_k(f_l))}, \quad \forall k=0,\ldots,n.$$

Most of the results on convergence properties of aggregation
methods are obtained for the regression and the gaussian white
noise models. Nevertheless, \citet{cat:97,catbook:04},
\citet{dl:01}, \citet{yang:00}, \citet{z:03} and \citet{rt:04}
have explored the performances of aggregation procedures in the
density estimation framework. Most of them have established upper
bounds for some procedure and do not deal with the problem of
optimality of their procedures. To our knowledge, lower bounds for
the performance of aggregation methods in density estimation are
available only in \citet{rt:04}. Their results are obtained with
respect to the mean squared risk. \citet{cat:97} and
\citet{yang:00} construct procedures and give convergence rates
w.r.t. the KL loss. One aim of this paper is to prove optimality
of one of these procedures w.r.t. the KL loss. Lower bounds w.r.t.
the Hellinger's distance and $L_1$-distance (stated in Section
\ref{sectionlowerbound}) and some results of \citet{bir:04} and
\citet{dl:01} (recalled in Section \ref{sectionupperbound})
suggest that the rates of convergence obtained in Theorem
\ref{theohellinger} and \ref{theoL1} are optimal in the sense
given in Definition \ref{definitionoptimality}. In fact, an
approximate bound can be achieved, if we allow the leading term in
the RHS of the oracle inequality (i.e. in the upper bound) to be
multiplied by a constant greater than one.

The paper is organized as follows. In Section
\ref{sectiondefinitionresult} we give a Definition of optimality,
for a rate of aggregation and for an aggregation procedure, and
our main results. Lower bounds, for different loss functions, are
given in Section \ref{sectionlowerbound}. In Section
\ref{sectionupperbound}, we recall a result of \citet{yang:00}
about an exact oracle inequality satisfied by the aggregation
procedure introduced in (\ref{esti}).

\par

\section{Main definition and main results}\label{sectiondefinitionresult}

To evaluate the accuracy of a density estimator we use the
Kullback-Leiber (KL) divergence, the Hellinger's distance and the
$L_1$-distance as loss functions. The {\emph{KL divergence}} is
defined for all densities $f$, $g$ w.r.t. a $\sigma-$finite
measure $\nu$ on a space $\cX$, by $$K(f|g)=\left\{
\begin{array}{ll}
\int_{\cal{X}} \log \left( \frac{f}{g}\right)f d \nu & \mbox {if } P_f \ll P_g ;\\
+\infty & \mbox{otherwise},
\end{array}\right.$$ where $P_f$ (respectively $P_g$) denotes the probability distribution of density $f$
(respectively $g$) w.r.t. $\nu$. {\emph{Hellinger's distance}} is
defined for all non-negative measurable functions $f$ and $g$ by
$$H(f,g)=\left\| \sqrt{f}-\sqrt{g} \right\|_{2},$$ where the
$L_2$-norm is defined by $
\|f\|_2=\left(\int_{\cal{X}}f^2(x)d\nu(x)\right)^{1/2}$ for all
functions $f \in L_2(\cal{X},\nu)$. The {\emph{$L_1$-distance}} is
defined for all measurable functions $f$ and $g$ by
$$v(f,g)=\int_{\cX}|f-g|d\nu.$$

The main goal of this paper is to find optimal rate of aggregation
in the sense of the definition given below. This definition is an
analog, for the density estimation problem, of the one in
\citet{tsy:03} for the regression problem.

\begin{defn}\label{definitionoptimality}
Take $M \geq 2$ an integer, $\cF$ a set of densities on
$(\cX,\cA,\nu)$ and $\cF_0$ a set of functions on $\cX$ with
values in $\mathbb{R}$ such that $\cF\subseteq\cF_0$. Let $d$ be a
loss function on the set $\cF_0$. A sequence of positive numbers $
(\psi _n (M))_{n \in \mathbb{N}^* }$ is called \textbf{optimal
rate of aggregation of M functions in $(\cF_0,\cF)$ w.r.t. the
loss $d$} if :
\begin{enumerate}[(i)]
\item There exists a constant $C< \infty $, depending only on
$\cF_0,\cF \mbox{ and } d $, such that for all functions
$f_1,\ldots,f_M$ in $\cF_0$ there exists an estimator $ \tilde
f_n$ (aggregate) of $f$  such that
\begin{equation}\label{upperbound}
\sup_{f \in \cF } \left[ \mathbb{E}_f \left[
d(f,\tilde{f_n})\right] - \min_{i=1, \ldots , M} d(f , f_i )
\right] \leq C \psi_n(M), \quad \forall n\in\mathbb{N}^{*}.
\end{equation}
\item There exist some functions $f_1,\ldots,f_M$ in $\cF_0$ and
$c>0$
 a constant independent of $M$ such that for all estimators $\hat f_n$
 of $f$,
\begin{equation}\label{lowerbound}
\sup_{f\in \cF } \left[ \mathbb{E}_f \left[ d(f ,
\hat{f_n})\right] - \min_{i=1, \ldots , M} d(f , f_i )  \right]
\geq c \psi_n(M), \quad \forall n\in\mathbb{N}^*.
\end{equation}
\end{enumerate}
Moreover, when the inequalities (\ref{upperbound}) and
(\ref{lowerbound}) are satisfied, we say that the procedure $
\tilde f_n$, appearing in (\ref{upperbound}), is an {\bf optimal
aggregation procedure w.r.t. the loss $d$.}
\end{defn}
Let $A>1$ be a given number. In this paper we are interested in
the estimation of densities lying in
\begin{equation}\label{set}\cF(A)= \left\{ \mbox{densities bounded by
} A \right\}\end{equation}and, depending on the used loss
function, we aggregate functions in $\cF_0$ which can be:
\begin{enumerate}
\item $\cF_K(A)= \left\{ \mbox{densities bounded by } A \right\}$
for KL divergence, \item $\cF_H(A)=\left\{ \mbox{non-negative
measurable functions bounded by } A \right\}$ for Hellinger's
distance, \item $\cF_v(A)=\left\{ \mbox{measurable functions
bounded by } A \right\}$ for the $L_1$-distance.
\end{enumerate}

The main result of this paper, obtained by using Theorem
\ref{theoyang} and assertion (\ref{ineq1}) of Theorem
\ref{theokullback}, is the following Theorem.
\begin{theo}\label{theooptimal} Let $A>1$. Let $M$ and $n$ be two integers such that
$\log M \leq 16(\min(1,A-1))^2 n$. The sequence
$$\psi_n(M)=\frac{\log M}{n}$$ is an optimal rate of aggregation
of $M$ functions in $(\cF_K(A),\cF(A))$ (introduced in
(\ref{set})) w.r.t. the KL divergence loss. Moreover, the
aggregation procedure with exponential weights, defined in
(\ref{esti}), achieves this rate. So, this procedure is an optimal
aggregation procedure w.r.t. the KL-loss.
\end{theo}

Moreover, observing Theorem \ref{theobirge} and the result of
\citet{dl:01} (recalled at the end of Section
\ref{sectionupperbound}), the rates obtained in Theorems
\ref{theohellinger} and \ref{theoL1}: $$\left( \frac{\log
M}{n}\right)^{\frac{q}{2}}$$ are near optimal rate of aggregation
for the Hellinger's distance and the $L_1$-distance to the power
$q$, where $q>0$, if we allow the leading term "$\min_{i=1, \ldots
, M} d(f , f_i )$"  to be multiplied by a constant greater than
one, in the upper bound and the lower bound.

\section{Lower bounds}\label{sectionlowerbound}
To prove lower bounds of type (\ref{lowerbound}) we use the
following lemma on minimax lower bounds which can be obtained by
combining Theorems 2.2 and 2.5 in \citet{tsybook:04}. We say that
$d$ is a {\emph semi-distance} on $\Theta$ if $d$ is symmetric,
satisfies the triangle inequality and $d(\theta,\theta)=0$.

\begin{lem}\label{lemlowerbound}
Let $d$ be a semi-distance  on the set of all densities on
$(\cX,\cA,\nu)$ and $w$ be a non-decreasing function defined on
$\mathbb{R}_{+}$ which is not identically $0$. Let
$(\psi_n)_{n\in\mathbb{N}}$ be a sequence of positive numbers. Let
$\cal{C}$ be a finite set of densities on $(\cX,\cA,\nu)$ such
that $card(\cC)=M \geq 2$,$$ \forall f\not= g \in \cC , \, d(f,g)
\geq 4 \psi_n>0 ,$$ and the KL divergences $K(P_f^{\otimes n}|
P_g^{\otimes n}) $, between the product probability measures
corresponding to densities $f$ and $g$ respectively, satisfy, for
some $f_0 \in \cC$,
$$ \forall f \in \cC,\, K(P_f^{\otimes n} | P_{f_0}^{\otimes n} )\leq (1/16) \log (M).$$
Then,
$$ \inf_{\hat{f_n}} \sup_{f \in \cal{C}} \mathbb{E}_f \left[ w(\psi
_n ^{-1}d( \hat{f_n},f )) \right]\geq c_1,$$ where
$\inf_{\hat{f_n}}$denotes the infimum over all estimators based on
a sample of size $n$ from an unknown distribution with density $f$
and $c_1>0$ is an absolute constant.
\end{lem}

Now, we give a lower bound of the form (\ref{lowerbound}) for the
three different loss functions introduced in the beginning of the
section. Lower bounds are given in the problem of estimation of a
density on $\mathbb{R}^d$, namely we have $\cX=\mathbb{R}^d$ and
$\nu$ is the Lebesgue measure on $\mathbb{R}^d$.

\begin{theo}\label{theohellinger}
Let $M$ be an integer greater than $2$, $A>1$ and $q>0$ two
numbers. We have for all integers $n$ such that $\log M \leq
16(\min{(1,A-1)})^2 n$,
$$ \sup_{f_1,\ldots,f_M\in\cF_H(A)} \inf_{\hat f_n}  \sup_{f\in \cF(A)}
\left[ \mathbb{E}_f \left[ H(\hat{f}_n,f)^q
\right]-\min_{j=1,\ldots,M} H(f_j,f)^q \right] \geq c
\left(\frac{\log M}{n}\right)^{q/2},$$ where $c$ is a positive
constant which depends only on $A$ and $q$. The sets $\cF(A)$ and
$\cF_H(A)$ are defined in (\ref{set}) when $\cX=\mathbb{R}^d$ and
the infimum is taken over all the estimators based on a sample of
size $n$.
\end{theo}

\noindent{\bf{Proof :}} For all densities  $f_1,\ldots,f_M$
bounded by $A$ we have,
$$\sup_{f_1,\ldots,f_M\in\cF_H(A)} \inf_{\hat f_n}  \sup_{f\in \cF(A)}
 \left[ \mathbb{E}_f \left[ H(\hat{f}_n,f)^q \right]-\min_{j=1,\ldots,M}H(f_j,f)^q \right]
 \geq \inf_{\hat f_n} \sup_{f\in \{f_1,\ldots,f_M \}}
 \mathbb{E}_f \left[ H(\hat{f}_n,f)^q \right].$$
Thus, to prove Theorem 1, it suffices to find $M$ appropriate
densities bounded by $A$ and to apply Lemma 1 with a suitable
rate.

We consider $D$ the smallest integer such that $2^{D/8}\geq M$ and
$\Delta=\{0,1\}^D$. We set $h_j(y)=h \left(y-(j-1)/D \right)$ for
all $y\in\mathbb{R}$, where $h(y)=(L/D)g(Dy)$ and
$g(y)=\1_{[0,1/2]}(y)-\1_{(1/2,1]}(y)$ for all $y\in\mathbb{R}$
and $L>0$ will be chosen later. We consider
$$ f_{\delta} (x)=\1_{[0,1]^d}(x)\left( 1+\sum_{j=1}^{D}\delta_j h_j(x_1)\right)
, \quad\forall x=(x_1,\ldots,x_d)\in\mathbb{R}^d,$$ for all
$\delta=(\delta_1,\ldots,\delta_D) \in\Delta$. We take $L$ such
that $L\leq D \min(1,A-1)$ thus, for all $\delta \in \Delta$,
$f_{\delta}$ is a density bounded by $A$. We choose our densities
$f_1,\ldots,f_M$ in
$\cB=\left\{f_{\delta}:\delta\in\Delta\right\}$, but we do not
take all of the densities of $\cB$ (because they are too close to
each other), but only a subset of $\cB$, indexed by a separated
set (this is a set where all the points are separated from each
other by a given distance) of $\Delta$ for the {\emph{Hamming
distance}} defined by
$\rho(\delta^1,\delta^2)=\sum_{i=1}^{D}I(\delta^1_i \not=
\delta^2_i)$ for all
$\delta^1=(\delta^1_1,\ldots\delta^1_D),\delta^2=(\delta^2_1,\ldots,\delta^2_D)
\in \Delta $. Since $\int_{\mathbb{R}}hd\lambda=0$, we have
\begin{eqnarray*}
H^2(f_{\delta^1},f_{\delta^2}) & = &
\sum_{j=1}^{D}\int_{\frac{j-1}{D}}^{\frac{j}{D}}I(\delta_j^1 \not=
\delta^2_j)
 \left(1-\sqrt{1+h_j(x)} \right)^2 dx \\
& = & 2\rho (\delta^1,\delta^2) \int_{0}^{1/D}\left(
1-\sqrt{1+h(x)}\right)dx ,
\end{eqnarray*}for all
$\delta^1=(\delta^1_1,\ldots,\delta^1_D),\delta^2=(\delta^2_1,\ldots,\delta^2_D)
\in \Delta$. On the other hand the function $\varphi (x) =1-\alpha
x^2 - \sqrt{1+x}$, where $\alpha = 8^{-3/2}$, is convex on $[-1,1]
$ and we have $ |h(x)|\leq L/D \leq 1$ so, according to Jensen,
$\int_{0}^{1}\varphi(h(x))dx\geq \varphi \left( \int_{0}^{1}h(x)dx
\right)$. Therefore $\int_0^{1/D} \left( 1-\sqrt{1+h(x)}\right)dx
\geq \alpha \int_0^{1/D}h^2(x)dx=(\alpha L^2)/D^3$, and we have
$$ H^2(f_{\delta^1},f_{\delta^2})\geq \frac{2\alpha L^2}{D^3} \rho(\delta^1,\delta^2),$$for all $
\delta^1,\delta^2 \in \Delta $. According to Varshamov-Gilbert,
cf. \citet[p. 89]{tsybook:04} or \citet{ih:80}, there exists a
$D/8$-separated set, called $N_{D/8}$, on $\Delta$ for the Hamming
distance such that its cardinal is higher than $2^{D/8}$ and
$(0,\ldots,0)\in N_{D/8}$. On the separated set $N_{D/8} $ we
have,
$$\forall \delta^1,\delta^2 \in N_{D/8} \, , \,  H^2(f_{\delta^1},f_{\delta^2})\geq \frac{\alpha L^2}{4D^2} .$$

In order to apply Lemma \ref{lemlowerbound}, we need to control
the KL divergences too. Since we have taken $N_{D/8}$ such that
$(0,\ldots,0) \in N_{D/8}$, we can control the KL divergences
w.r.t. $P_0 $, the Lebesgue measure on $[0,1]^d$. We denote by
$P_{\delta}$ the probability of density $f_\delta$ w.r.t. the
Lebesgue's measure on $\mathbb{R}^d$, for all $\delta \in \Delta$.
We have,
\begin{eqnarray*}
K(P_{\delta}^{\otimes n}|P_0^{\otimes n}) & = & n \int_{[0,1]^d} \log \left( f_{\delta}(x)\right)f_{\delta}(x) dx \\
& = & n \sum_{j=1}^{D}\int_{\frac{j-1}{D}}^{j/D}\log \left(
1+\delta_j h_j(x)\right)\left( 1+\delta_j h_j(x)\right)dx \\
& = & n
\left(\sum_{j=1}^{D}\delta_j\right)\int_{0}^{1/D}\log(1+h(x))(1+h(x))dx,
\end{eqnarray*}
for all $\delta=(\delta_1,\ldots,\delta_D) \in N_{D/8}$. Since
$\forall u
>-1 , \log(1+u) \leq u$, we have,
\begin{equation*}
K(P_{\delta}^{\otimes n}|P_0^{\otimes n})\leq n \left(
\sum_{j=1}^{D}\delta _j \right) \int_{0}^{1/D}(1+h(x))h(x)dx  \leq
nD\int_{0}^{1/D}h^{2}(x)dx=\frac{n L^2}{D^2}.
\end{equation*}

Since $\log M \leq 16(\min{(1,A-1)})^2 n$, we can take $L$ such
that $(nL^2)/D^2=\log(M)/16$ and still having $L \leq
D\min(1,A-1)$. Thus, for $L=(D/4) \sqrt{\log(M)/n}$, we have for
all elements $\delta^1,\delta^2$ in $N_{D/8}$,
$H^2(f_{\delta^1},f_{\delta^2})\geq (\alpha/64) (\log(M)/n)$ and
$\forall \delta \in N_{D/8}\, , \, K(P_\delta^{\otimes n} |
P_0^{\otimes n} ) \leq (1/16) \log(M).$

Applying Lemma 1 when $d$ is $H$, the Hellinger's distance, with
$M$ densities $f_1,\ldots,f_M$ in $\left\{f_{\delta} : \delta \in
N_{D/8} \right\}$ where $f_1=\1_{[0,1]^d}$ and the increasing
function $w(u)=u^q$, we get the  result.
\begin{flushright} $\blacksquare$ \end{flushright}

\begin{remark}
The construction of the family of densities $\left\{f_{\delta} :
\delta \in N_{D/8} \right\}$ is in the same spirit as the lower
bound of \citet{tsy:03}, \citet{rt:04} but, as compared to
\citet{rt:04},  we consider a different problem (model selection
aggregation) and as compared to \citet{tsy:03}, we study a
different model (density estimation). Also, our risk function is
different from those considered in these papers.
\end{remark}

Now, we give a lower bound for KL divergence. We have the same
residual as for square of Hellinger's distance.
\begin{theo}\label{theokullback}Let $M\geq2$ be an integer, $A>1$ and $q>0$.
We have, for any integer $n$ such that $\log M \leq
16(\min(1,A-1))^2 n$,
\begin{equation}\label{ineq1}
\sup_{f_1,\ldots,f_M \in\cF_K(A)}\inf_{\hat f_n} \sup_{
f\in\cF(A)}\left[ \mathbb{E}_f \left[ (K(f|\hat{f_n}))^q \right]
 -\min_{j=1,\ldots,M} (K(f|f_j))^q \right]
\geq c \left(\frac{\log M}{n}\right)^q, \end{equation}and
\begin{equation}\label{ineq2}
\sup_{f_1,\ldots,f_M \in\cF_K(A)}\inf_{\hat f_n}\sup_{ f\in\cF(A)}
\left[ \mathbb{E}_f \left[ (K(\hat{f_n}|f))^q \right]
-\min_{j=1,\ldots,M} (K(f_j|f))^q \right]\geq c \left(\frac{\log
M}{n}\right)^q,\end{equation}where $c$ is a positive constant
which depends only on $A$. The sets $\cF(A)$ and $\cF_K(A)$ are
defined in (\ref{set}) for $\cX=\mathbb{R}^d$.
\end{theo}

\noindent{\bf{Proof :}} Proof of the inequality (\ref{ineq2}) of
Theorem \ref{theokullback} is similar to the one for
(\ref{ineq1}). Since we have for all densities $f$ and $g$,
$$K(f|g)\geq H^2(f,g),$$  \citep[a proof is given in][p. 73]{tsybook:04}, it suffices to note that, if
$f_1,\ldots,f_M$ are densities bounded by $A$ then,
$$\sup_{f_1,\ldots,f_M \in\cF_K(A)}\inf_{\hat f_n} \sup_{ f
\in\cF(A)} \left[ \mathbb{E}_f \left[ (K(f|\hat{f_n}))^q \right]
-\min_{j=1,\ldots,M} (K(f|f_i))^q \right]$$ $$ \geq \inf_{\hat
f_n} \sup_{f\in \{f_1,\ldots,f_M \}} \left[ \mathbb{E}_f \left[
(K(f|\hat{f_n}))^q \right] \right] \geq  \inf_{\hat f_n}
\sup_{f\in \{f_1,\ldots,f_M \}} \left[ \mathbb{E}_f \left[
H^{2q}(f,\hat{f_n}) \right] \right],$$to get the result by
applying Theorem \ref{theohellinger}.
\begin{flushright} $\blacksquare$ \end{flushright}

With the same method as Theorem 1, we get the result below for the
$L_1$-distance.

\begin{theo}\label{theoL1} Let $M\geq2$ be an integer, $A
> 1$ and $q>0$. We have for any integers $n$ such that $\log M
\leq 16(\min(1,A-1))^2 n$,
$$ \sup_{f_1,\ldots,f_M \in\cF_v(A)}\inf_{\hat f_n} \sup_{f \in\cF(A)} \left[ \mathbb{E}_f \left[ v(f,\hat{f_n})^q \right]
 -\min_{j=1,\ldots,M} v(f,f_i)^q \right]\geq c \left(\frac{\log M}{n}\right)^{q/2}$$
where $c$ is a positive constant which depends only on $A$. The
sets $\cF(A)$ and $\cF_v(A)$ are defined in (\ref{set}) for
$\cX=\mathbb{R}^d$.
\end{theo}

\noindent{\bf{Proof :}} The only difference with Theorem
\ref{theohellinger} is in the control of the distances. With the
same notations as the proof of Theorem \ref{theohellinger}, we
have,
$$v(f_{\delta^1},f_{\delta^2})=\int_{[0,1]^d}
|f_{\delta^1}(x)-f_{\delta^2}(x)|dx=
\rho(\delta^1,\delta^2)\int_0^{1/D} |h(x)|dx
=\frac{L}{D^2}\rho(\delta^1,\delta^2),$$for all $\delta^1,\delta^2
\in\Delta$. Thus, for $L=(D/4)\sqrt{\log(M)/n}$ and $N_{D/8}$, the
$D/8$-separated set of $\Delta$ introduced in the proof of Theorem
\ref{theohellinger}, we have,
$$v(f_{\delta^1},f_{\delta^2}) \geq
\frac{1}{32}\sqrt{\frac{\log(M)}{n}}, \quad\forall \delta^1,
\delta^2 \in N_{D/8} \mbox{ and } K(P_{\delta}^{\otimes n}
|P_0^{\otimes n})\leq \frac{1}{16} \log(M), \quad\forall \delta
\in \Delta.$$ Therefore, by applying Lemma 1 to the $L_1$-distance
with $M$ densities $f_1,\ldots,f_M$ in $\left\{f_{\delta} : \delta
\in N_{D/8} \right\}$ where $f_1=\1_{[0,1]^d}$ and the increasing
function $w(u)=u^q$, we get the result.
\begin{flushright} $\blacksquare$ \end{flushright}

\par

\section{Upper bounds}\label{sectionupperbound}

In this section we use an argument in \citet{yang:00} \citep[see
also][]{catbook:04} to show that the rate of the lower bound of
Theorem \ref{theokullback} is an optimal rate of aggregation with
respect to the KL loss. We use an aggregate constructed by Yang
(defined in (\ref{esti})) to attain this rate. An upper bound of
the type (\ref{upperbound}) is stated in the following Theorem.
Remark that Theorem \ref{theoyang} holds in a general framework of
a measurable space $(\cX,\cA)$ endowed with a $\sigma$-finite
measure $\nu$.

\begin{theo}[Yang]\label{theoyang}
Let $X_1,\ldots,X_n$ be $n$ observations of a probability measure
on $(\cX,\cA)$ of density $f$ with respect to $\nu$. Let
$f_1,\ldots,f_M$ be $M$ densities on $(\cal{X},\cA,\nu)$. The
aggregate $\tilde{f}_n$, introduced in (\ref{esti}), satisfies,
for any underlying density $f$,
\begin{equation}\label{upperboundyang}\mathbb{E}_f \left[ K(f|\tilde{f}_n)
\right] \leq \min_{j=1,\ldots,M}K(f|f_j) +
\frac{\log(M)}{n+1}.\end{equation}
\end{theo}

\noindent{\bf{Proof :}} Proof follows the line of \citet{yang:00},
although he does not state the result in the form
(\ref{upperbound}), for convenience we reproduce the argument
here. We define $
\hat{f}_k(x;X^{(k)})=\sum_{j=1}^{M}w_j^{(k)}f_j(x), \ \forall
k=1,\ldots,n$ (where $w_j^{(k)}$ is defined in (\ref{weights}) and
$x^{(k)}=(x_1,\ldots,x_k)$ for all $k\in\mathbb{N}$ and
$x_1,\ldots,x_k\in\cX$) and
 $\hat{f}_0(x;X^{(0)})=(1/M)\sum_{j=1}^M f_j(x)$ for all
$x\in\cX$. Thus, we have \begin{equation*}
\tilde{f}_n(x;X^{(n)})=\frac{1}{n+1}\sum_{k=0}^n
\hat{f}_k(x;X^{(k)}).\end{equation*} Let $f$ be a density on
$(\cX,\cA,\nu)$. We have
\begin{eqnarray*}
\sum_{k=0}^{n} \mathbb{E}_f \left[ K(f|\hat{f}_k ) \right] & = & \sum_{k=0}^{n} \int_{\cX^{k+1}} \log \left(\frac{f(x_{k+1})}{\hat{f}_k(x_{k+1};x^{(k)})} \right)\prod_{i=1}^{k+1}f(x_i)d\nu^{\otimes (k+1)}(x_1,\ldots,x_{k+1})\\
& = & \int_{\cX^{n+1}} \left( \sum_{k=0}^{n}\log \left(\frac{f(x_{k+1})}{\hat{f}_k(x_{k+1};x^{(k)})} \right) \right)\prod_{i=1}^{n+1}f(x_i)d\nu^{\otimes (n+1)}(x_1,\ldots,x_{n+1})\\
& = &\int_{\cX^{n+1}}\log \left( \frac{f(x_1)\ldots f(x_{n+1})}
{\prod_{k=0}^{n}\hat{f}_k(x_{k+1};x^{(k)})}\right)\prod_{i=1}^{n+1}f(x_i)d\nu^{\otimes
(n+1)}(x_1,\ldots,x_{n+1}),
\end{eqnarray*}
but $\prod_{k=0}^{n}\hat{f}_k(x_{k+1};x^{(k)})= (1/M)
\sum_{j=1}^{M}f_j(x_1)\ldots f_j(x_{n+1}),\forall
x_1,\ldots,x_{n+1}\in\cX$ thus,
$$\sum_{k=0}^{n} \mathbb{E}_f \left[ K(f|\hat{f}_k ) \right]=\int_{\cX^{n+1}}\log \left( \frac{f(x_1)\ldots f(x_{n+1})}{\frac{1}{M} \sum_{j=1}^{M}f_j(x_1)\ldots f_j(x_{n+1})}\right)\prod_{i=1}^{n+1}f(x_i)d\nu^{\otimes (n+1)}(x_1,\ldots,x_{n+1}),$$
moreover $x \longmapsto \log(1/x)$ is a decreasing function so,
\begin{equation*}
\sum_{k=0}^{n} \mathbb{E}_f \left[ K(f|\hat{f}_k ) \right] \leq
\min_{j=1,\ldots,M} \left\{ \int_{\cX^{n+1}}\log \left(
\frac{f(x_1)\ldots f(x_{n+1})}{\frac{1}{M} f_j(x_1)\ldots
f_j(x_{n+1})}\right) \prod_{i=1}^{n+1}f(x_i)d\nu^{\otimes
(n+1)}(x_1,\ldots,x_{n+1})\right\} \end{equation*}
\begin{equation*}
\leq \log M + \min_{j=1,\ldots, M} \left\{ \int_{\cX^{n+1}}\log
\left( \frac{f(x_1)\ldots f(x_{n+1})}{f_j(x_1)\ldots
f_j(x_{n+1})}\right)\prod_{i=1}^{n+1}f(x_i)d\nu^{\otimes
(n+1)}(x_1,\ldots,x_{n+1})\right\},
\end{equation*} finally we have,
\begin{equation}\label{demoyang1}
\sum_{k=0}^{n} \mathbb{E}_f \left[ K(f|\hat{f}_k ) \right] \leq
\log M +(n+1) \inf_{j=1,\ldots, M } K(f|f_j).
\end{equation}
On the other hand we have,
$$\mathbb{E}_f \left[ K(f| \tilde{f}_{n}) \right]=
 \int_{\cX ^{n+1}} \log \left( \frac{f(x_{n+1})}{\frac{1}{n+1} \sum_{k=0}^{n}\hat{f}_{k} (x_{n+1};x^{(k)})}\right)
  \prod_{i=1}^{n+1}f(x_i)d\nu^{\otimes (n+1)}(x_1,\ldots,x_{n+1}),$$
and $x \longmapsto \log(1/x)$ is convex, thus,
\begin{equation}\label{demoyang2}
\mathbb{E}_f \left[ K(f| \tilde{f}_{n}) \right] \leq \frac{1}{n+1}
\sum_{k=0}^{n} \mathbb{E}_f \left[ K(f|\hat{f}_k ) \right].
\end{equation}
Theorem \ref{theoyang} follows by combining (\ref{demoyang1}) and
(\ref{demoyang2}).
\begin{flushright} $\blacksquare$ \end{flushright}

Birg\'e constructs estimators, called \emph{T-estimators} (the
''T'' is for ''test''), which are adaptive in aggregation
selection model of $M$ estimators with a residual proportional at
$\left(\log M/n \right)^{q/2}$ when Hellinger and $L_1$-distances
are used to evaluate the quality of estimation (cf.
\citet{bir:04}). But it does not give an optimal result as Yang,
because there is a constant greater than 1 in front of the main
term $\min_{i=1,\ldots, M}d^q(f,f_i)$ where $d$ is the Hellinger
distance or the $L_1$ distance. Nevertheless, observing the proof
of Theorem \ref{theohellinger} and \ref{theoL1}, we can obtain
$$\sup_{f_1,\ldots,f_M\in\cF(A)}\inf_{\hat{f}_n}\sup_{f\in \cF(A) } \left[ \mathbb{E}_f \left[ d(f ,
\hat{f_n})^q\right] - C(q)\min_{i=1, \ldots , M} d(f , f_i )^q
\right] \geq c\left(\frac{\log M}{n}\right)^{q/2},$$ where $d$ is
the Hellinger or $L_1$-distance, $q>0$ and $A>1$. The constant
$C(q)$ can be chosen equal to the one appearing in the following
Theorem. The same residual appears in this lower bound and in the
upper bounds of Theorem \ref{theobirge}, so we can say that
$$\left(\frac{\log M}{n}\right)^{q/2}$$is near optimal rate of aggregation w.r.t. the
Hellinger distance or the $L_1$-distance to the power $q$. We
recall Birg\'e's results in the following Theorem.

\begin{theo}[Birg\'e]\label{theobirge}
If we have $n$ observations of a probability measure of density
$f$ w.r.t. $\nu$ and $f_1,\ldots,f_M$ densities on
$(\cX,\cA,\nu)$, then there exists an estimator $\tilde{f}_n$ (
T-estimator) such that for any underlying density $f$ and $q>0$,
we have
$$\mathbb{E}_f \left[H(f,\tilde{f}_n)^q\right]\leq C(q)
\left( \min_{j=1,\ldots, M} H(f,f_j)^q+\left(\frac{\log
M}{n}\right)^{q/2} \right),$$ and for the $L_1$-distance we can
construct an estimator $\tilde{f}_n$ which satisfies :
$$\mathbb{E}_f \left[v(f,\tilde{f}_n)^q\right]\leq C(q)
\left( \min_{j=1,\ldots, M} v(f,f_j)^q+\left(\frac{\log
M}{n}\right)^{q/2} \right),$$ where $C(q)>0$ is a constant
depending only on $q$.
\end{theo}

An other result, which can be found in \citet{dl:01}, states that
the minimum distance estimate proposed by Yatracos (1985) (cf.
\citet[p.~59]{dl:01}) achieves the same aggregation rate as in
Theorem \ref{theobirge} for the $L_1$-distance with $q=1$. Namely,
for all $f,f_1,\ldots,f_M\in\cF(A)$,
$$\mathbb{E}_f \left[v(f,\breve{f}_n)\right]\leq
3\min_{j=1,\ldots, M} v(f,f_j)+\sqrt{\frac{\log M}{n}},$$ where
$\breve{f}_n$ is the estimator of Yatracos defined by
$$\breve{f}_n={\rm arg}\min_{f\in\{f_1,\ldots,f_M\}} \sup_{A\in\cA}
\left\vert \int_A f-\frac{1}{n}\sum_{i=1}^n\1_{\{X_i\in A\}}
\right\vert,$$ and $\cA=\left\{\{x:f_i(x)>f_j(x)\}:1\leq i, j \leq
M \right\}.$


\end{document}